\documentclass[10pt,reqno]{amsart}
\usepackage{amssymb,mathrsfs,graphicx,extpfeil,amsmath,bbm}
\usepackage{ifthen,latexsym,float,colortbl}
\usepackage[margin=0.7in]{geometry}
\usepackage{caption}
\usepackage{dsfont}
\usepackage{arydshln}
\usepackage{sidecap}
\usepackage{rotating,dsfont,stackengine}
\usepackage{hyperref}
\usepackage{refcheck}
\usepackage{enumerate}

\usepackage{colortbl}
\definecolor{black}{rgb}{0.0, 0.0, 0.0}
\definecolor{red}{rgb}{1.0, 0.5, 0.5}
\provideboolean{shownotes} 
\setboolean{shownotes}{true} 
\newcommand{\margnote}[1]{
	\ifthenelse{\boolean{shownotes}}%
	{\marginpar{\raggedright\tiny\texttt{#1}}}%
	{}%
}
\newcommand{\hole}[1]{
	\ifthenelse{\boolean{shownotes}}%
	{\begin{center} \fbox{ \rule {.25cm}{0cm} \rule[-.1cm]{0cm}{.4cm}
				\parbox{.85\textwidth}{\begin{center} \texttt{#1}\end{center}} \rule
				{.25cm}{0cm}}\end{center}} {} }

\graphicspath{{pics/}}

\title[Equilibrium parameters of the Marle model for polyatomic gases]
{Determination of equilibrium parameters of the Marle model for  polyatomic gases}

\author[Hwang]{Byung-Hoon Hwang}
\address[Byung-Hoon Hwang]{\newline Department of Mathematics Education\newline
	Sangmyung University, 20 Hongjimun 2-gil, Jongno-Gu, Seoul 03016, Republic of Korea}
\email{bhhwang@smu.ac.kr}

\numberwithin{equation}{section}

\newtheorem{theorem}{Theorem}[section]

\newtheorem{remark}{Remark}[section]

\newcommand{\R}{\mathbb R}

\newcommand{\I}{\mathcal I}

\newcommand{\bq}{\begin{equation}}
	\newcommand{\eq}{\end{equation}}

\newcommand{\pa}{\partial}

\newcommand{\intr}{\int_{\R^3}}
\newcommand{\inti}{\int_0^\infty}

\makeatletter
\def\moverlay{\mathpalette\mov@rlay}
\def\mov@rlay#1#2{\leavevmode\vtop{%
		\baselineskip\z@skip \lineskiplimit-\maxdimen
		\ialign{\hfil$\m@th#1##$\hfil\cr#2\crcr}}}
\newcommand{\charfusion}[3][\mathord]{
	#1{\ifx#1\mathop\vphantom{#2}\fi
		\mathpalette\mov@rlay{#2\cr#3}
	}
	\ifx#1\mathop\expandafter\displaylimits\fi}
\makeatother

\makeatletter
\newcommand*{\rom}[1]{\expandafter\@slowromancap\romannumeral #1@}
\makeatother
\makeatletter
\@namedef{subjclassname@2020}{\textup{2020} Mathematics Subject Classification}
\makeatother
\begin{document}
	\allowdisplaybreaks
	
	\date{\today}
	
	\subjclass[2020]{}

	\begin{abstract} 
			 The BGK model is a relaxation-time approximation of the celebrated Boltzmann equation, and the Marle model is a direct extension of the BGK model in a relativistic framework. In this paper, we introduce the Marle model for polyatomic gases based on the J\"{u}ttner distribution devised in [Ann. Phys., 377, (2017), 414--445], and show the existence of a unique set of equilibrium parameters of the J\"{u}ttner distribution.
	\end{abstract}
	
	\maketitle \centerline{\date}

	%
	%
	%
	%
	\setcounter{equation}{0}
	\section{Introduction}
\subsection{Relativistic extended thermodynamics of polyatomic molecules}
Recently, in 2017, a relativistic extended thermodynamics (RET) of rarefied polyatomic gases was discussed in \cite{PR17}, where the relativistic Maxwellian, also called the J\"{u}ttner distribution, was derived for polyatomic gases for the first time. To be precise, let $f\equiv f(x^\alpha , p^\beta,\mathcal{I} )$ be the momentum distribution function of relativistic particles on the phase space point $(x^\alpha,p^\beta)$ with the internal energy $\mathcal{I}\ge 0$ due to the rotation and vibrations of particles, where  $x^\alpha=(ct,x)\in\mathbb{R}^+\times \mathbb{R}^3$ is the space-time coordinate, and $p^\beta=(\sqrt{(mc)^2+|p|},p)\in\mathbb{R}^+\times \mathbb{R}^3$ is the four-momentum. Here $c$ is the speed of light, $m$ is the rest mass of a particle, and Greek indices
run from $0$ to $3$. Throughout this paper, the metric tensor $g_{\alpha\beta}$ and its inverse $g^{\alpha\beta}$  are given by 
$$
g_{\alpha\beta}=g^{\alpha\beta}=\text{diag}(1,-1,-1,-1)
$$
and we use the raising and lowering indices as
$$
g_{\alpha\mu}p^\mu=p_\alpha,\qquad g^{\alpha\mu}p_\mu=p^\alpha,
$$
which implies $p_\alpha=(p^0,-p)$. Then it follows from the Einstein summation convention that
$$
p^\mu q_\mu=p_\mu q^\mu=p^0q^0-\sum_{i=1}^3 p^iq^i.
$$
For $f\equiv f(t,x,p,\I)$, macroscopic descriptions are given by the particle-particle flux $V^\mu$ and energy-momentum tensor $T^{\mu\nu}$ \cite{PR17}:
\begin{equation}\label{VT}
V^\mu=mc\int_{\mathbb{R}^3}\int_0^\infty p^\mu f \phi(\mathcal{I}) \,d\mathcal{I}\,\frac{dp}{p^0},\qquad T^{\mu\nu}=\frac{1}{mc} \int_{\mathbb{R}^3}\int_0^\infty p^\mu p^\nu f\left( mc^2 + \mathcal{I} \right)   \phi(\mathcal{I})\,d\mathcal{I}\,\frac{dp}{p^0},
\end{equation}
where $\phi(\mathcal{I})\ge 0$ is the state density so that  $\phi(\mathcal{I}) \, d  \mathcal{I}$ represents the number of the internal states of a molecule having the internal energy between $\mathcal{I}$ and $\mathcal{I}+d \mathcal{I}$.  The form of $\phi(\mathcal{I})$ may be determined differently depending on the physical context, but it should be able to recover the state density for classical particles in the non-relativistic limit:$$
	 \lim_{c\to\infty}\phi(\I)=I^\sigma\qquad \mbox{with}\qquad \sigma=\frac{f^i-2}{2},
	 $$
where $I$ denotes the variable representing the internal energy for classical particles and  $f^i\ge 0$ is the internal degrees of freedom due to the internal motion of molecules (see \cite[Sec 4.1]{PR17} for details). For monatomic gases, $\sigma=-1$. In this paper, we choose $\phi(\mathcal{I})$ of the form 
$$
\phi(\mathcal{I})=\mathcal{I}^{\sigma}\ (\sigma >-1),
$$
 which is employed in \cite{PR18} to propose another type of relativistic BGK model for polyatomic gases. The entropy four-vector $h^\mu$ is given by
$$
h^\mu=-k_B c\intr \inti p^\mu f\ln f\phi(\I)\,d\I \frac{dp}{p^0},
$$
where $k_B$ denotes the Boltzmann constant. Finally, we introduce the J\"{u}ttner distribution for polyatomic gases:
\begin{equation*} 
f_E\equiv f_E(n,U^\mu,\gamma)=\frac{n}{M(\gamma)}e^{-\left(1+\frac{\mathcal{I}}{mc^2}\right)\frac{\gamma}{mc^2} U^\mu p_\mu }\qquad \mbox{with}\qquad M(\gamma) =\int_{\mathbb{R}^3}\int_0^\infty  e^{-\left(1+\frac{\mathcal{I}}{mc^2}\right)\frac{\gamma}{mc} p^0 }  \phi(\mathcal{I}) \,d\mathcal{I}\,dp,
\end{equation*}
which maximizes the entropy $h:=h^\alpha U_\alpha$ (for details, see \cite[Sec.4]{PR17}). Here the equilibrium parameters $n,U^\mu=(\sqrt{c^2+|U|^2},U),$ and $\gamma:=\frac{mc^2}{k_BT}$ are unknown functions of $t$ and $x$ representing the equilibrium density, four-velocity, and the ratio between the rest energy of a particle and  the product of the Boltzmann constant and equilibrium temperature $T$ respectively.

\subsection{Marle model for polyatomic gases} The BGK model \cite{BGK54} is a relaxation-time approximation of the celebrated Boltzmann equation, which has been widely used in physics and engineering for practical purposes. For monatomic molecules, a direct extension of the BGK model in the relativistic framework was first proposed by Marle  \cite{Marle65, Marle69} based on the Eckart frame \cite{Eckart40}. In this paper, we present the direct application of RET \cite{PR17} to the Marle model \cite{Marle65, Marle69}, which reads
\begin{align}\label{PMarle}
\pa_t f+\hat{p}\cdot\nabla_x f=\frac{cm}{\tau (1+\frac{\I}{mc^2})p^0}(f_E-f) =:Q(f)
\end{align}
where $\hat{p}:=cp/p^0$ is the normalized momentum and $\tau$ denotes the relaxation time in the rest frame. In the Eckart frame \cite{Eckart40}, the particle-particle flux $V_f^\mu$ of \eqref{VT} can be expressed by the observable quantities as $\displaystyle 
V_f^\mu=mn_f U_f^\mu$ where $mn_f$ denotes the number density and $U_f^\mu=(\sqrt{c^2+|U_f|^2},U_f)$ the Eckart four-velocity defined by
\begin{align}\label{eckart}
\begin{split}
n_f^2&=\left(\int_{\mathbb{R}^3}\int_0^\infty  f \phi(\mathcal{I}) \,d\mathcal{I}\,dp\right)^2-\sum_{i=1}^3\left(\int_{\mathbb{R}^3}\int_0^\infty p^i f \phi(\mathcal{I}) \,d\mathcal{I}\,\frac{dp}{p^0}\right)^2,\cr 
U_f^\mu &=\frac {c}{n_f}\int_{\mathbb{R}^3}\int_0^\infty p^\mu f \phi(\mathcal{I}) \,d\mathcal{I}\,\frac{dp}{p^0}
\end{split}\end{align}
respectively. Since $U^\mu_f$ has a constant length, i.e. $U^\mu_f U_{f\mu}=c^2$, the number density $n_f$ can be written with respect to $U_f^\mu $ as
\begin{equation}\label{n}
n_f= \frac1c \int_{\mathbb{R}^3}\int_0^\infty U_{f\mu} p^\mu f \phi(\mathcal{I}) \,d\mathcal{I}\,\frac{dp}{p^0}.
\end{equation}
Note that  the head term $(1+\frac{\I}{mc^2})^{-1}$ on the right-hand side of \eqref{PMarle} is considered for consistency with the particle-particle flux $V_f^\mu$  in the application of the Eckart frame (see \eqref{2}).
\subsection{Main result} In relativistic BGK-type models, the equilibrium parameters $n,U^\mu$ and $\gamma$ of the J\"{u}ttner distribution are determined in a way that the equation satisfies the conservation laws of $V^\mu$ and $T^{\mu\nu}$. In this process, $\gamma$ is often defined through a nonlinear relation due to the relativistic nature. Therefore, in order to study the equations rigorously,  it is necessary to show that $\gamma$ can be uniquely determined by the distribution $f$, see \cite{BCNS12, Hwang22, HLY23, HRY22} for similar problems.  We also refer to \cite{Calvo12, SS11} for the hydrodynamic limit of the relativistic Boltzmann equation regarding the range of $\gamma$.  The aim of this paper is to prove that the equilibrium parameters of $f_E$ in the polyatomic Marle model \eqref{PMarle} can be uniquely determined in a way that \eqref{PMarle} satisfies the fundamental properties of the Boltzmann equation, the conservation laws and $H$-theorem. Our main result is as follows.
\begin{theorem}\label{main}
	Let $f$ be non-negative and not trivially zero so that $V_f^\mu$ exists.  Assume that the state density is chosen as $\phi(\I)=\I^\sigma$ with $\sigma >-1$. Then there exists a unique set of equilibrium parameters $n,U^\mu$ and $\gamma$ of $f_E$ satisfying the following identities:
	\begin{equation}\label{cancellation}
\intr\inti Q(f)\phi(\I)\,d\I dp=0,\qquad \intr\inti p^\mu (mc^2+\I)Q(f)\phi(\I)\,d\I dp=0.
	\end{equation}
	Indeed, $n=n_f$, $U^\mu=U_f^\mu$, and $\gamma$ is determined by the nonlinear relation:
\begin{align}\label{gamma}
\frac{\widetilde{M}(\gamma)}{M(\gamma)}&=\frac {1}{n_f}\intr\inti  f\Big(1+\frac{\I}{mc^2}\Big)^{-1}\phi(\mathcal{I}) \,d\mathcal{I} \frac{dp}{p^0}.
\end{align}
where 
$$
\widetilde{M}(\gamma)=\intr\inti  e^{-\left(1+\frac{\mathcal{I}}{mc^2}\right)\frac{\gamma}{mc}p^0 }\Big(1+\frac{\I}{mc^2}\Big)^{-1}\phi(\mathcal{I}) \,d\mathcal{I} \frac{dp}{p^0}.
$$ 	 
\end{theorem}
\begin{remark}
Due to \eqref{cancellation}, the distribution function $f$ verifying \eqref{PMarle} satisfies the conservation laws of $V^\mu_f$ and $T_f^{\mu\nu}$:
\begin{align*}
\pa_\mu V_f^\mu &=m\int_{\mathbb{R}^3}\int_0^\infty \{\pa_t f+ \hat{p}\cdot\nabla_x f \}  \phi(\mathcal{I}) \,d\mathcal{I}\,dp=0,\cr 
 \pa_\nu T_f^{\mu\nu}&=\frac{1}{mc^2} \int_{\mathbb{R}^3}\int_0^\infty p^\mu \{\pa_t f+ \hat{p}\cdot\nabla_x f \}\left( mc^2 + \mathcal{I} \right)\phi(\mathcal{I})\,d\mathcal{I}\,dp=0,
\end{align*}
where $\pa_\mu:=\pa/\pa x^\mu$. Furthermore,  \eqref{cancellation} implies 
\begin{align*}
\intr\int_0^\infty Q(f)\ln f_E \phi(\I) \,d\I dp&=\ln \Big\{\frac{n}{M(\gamma)}\Big\}\intr\int_0^\infty Q(f) \phi(\I) \,d\I dp-\frac{\gamma}{m^2c^4}U_\mu \intr\int_0^\infty p^\mu (mc^2+\I)Q(f) \phi(\I) \,d\I dp\cr 
&=0,
\end{align*}
which gives
\begin{align*}
\intr\int_0^\infty Q(f)\ln f \phi(\I) \,d\I dp&=\intr\int_0^\infty Q(f)(\ln f -\ln f_E) \phi(\I) \,d\I dp\cr 
&=\frac{cm}{\tau }\intr\int_0^\infty (f_E-f)(\ln f -\ln f_E) \Big(1+\frac{\I}{mc^2}\Big)^{-1}\phi(\I) \,d\I \frac{dp}{p^0}\cr 
&=\frac{cm}{\tau }\intr\int_0^\infty f_E\Big(1-\frac{f}{f_E}\Big)\ln \Big(\frac{f}{f_E} \Big) \Big(1+\frac{\I}{mc^2}\Big)^{-1}\phi(\I) \,d\I \frac{dp}{p^0}.
\end{align*}
Using the following inequality:
$$
(1-x)\ln x\le 0\quad \mbox{for all } x> 0,
$$ 
one can obtain the $H$-theorem:  

$$
\pa_\mu h^\mu:=  -k_B \intr \inti \{ \pa_t+\hat{p}\cdot\nabla_x\} f\ln f\phi(\I)\,d\I dp \ge 0.
$$
\end{remark}
\begin{remark}  In this paper, the state density is chosen as $\phi(\I)=\I^\sigma$ $(\sigma>-1)$, which is used in the proof of Theorem \ref{main}, see \eqref{phi lower}. However, to obtain \eqref{phi lower}, other types of state densities also seem available if the following conditions are satisfied: (1) $\phi(\I)$ does not grow exponentially, (2) $\I\phi(\I)=0$ when $\I=0$, and (3) $\phi(\I)+\I\phi^\prime(\I)$ is governed by $(1+\frac{\I}{mc^2})\phi(\I)$ up to constant multiplication.  
\end{remark}
\smallskip
\section{Proof of theorem \ref{main}}
The proof is divided into two steps. We first show how the equilibrium parameters are determined by using the identities \eqref{cancellation} as constraints. And then we prove that $\gamma$ can be uniquely determined through the nonlinear relation \eqref{gamma}.
\noindent\newline
$\bullet $ Choice of $n,U^\mu,\gamma$: By definition of $Q(f)$, \eqref{cancellation} is written as
\begin{align}\label{1}
\intr \inti f_E\Big(1+\frac{\mathcal{I}}{mc^2}\Big)^{-1} \phi(\I)\, d\I \frac{dp}{p^0}&=	\intr \inti f\Big(1+\frac{\mathcal{I}}{mc^2}\Big)^{-1}\phi(\I)\, d\I \frac{dp}{p^0},
\end{align}
and
\begin{align}\label{2}
\intr \inti p^\mu f_E \phi(\I)\, d\I \frac{dp}{p^0}&=	\intr \inti p^\mu f\phi(\I)\, d\I \frac{dp}{p^0}.	
\end{align}
To simplify the integral of $f_E$, we introduce the Lorentz transformation	$\Lambda$: 
\begin{align*}
\Lambda=
\begin{bmatrix}
c^{-1}U^0 & -c^{-1}U^1 & -c^{-1}U^2 & -c^{-1}U^3 \cr
-U^1&  1+(U^0-1)\frac{(U^1)^2}{|U|^2}&(U^0-1)\frac{U^1U^2}{|U|^2}  &(U^0-1)\frac{U^1U^3}{|U|^2}  \cr
-U^2& (U^0-1)\frac{U^1U^2}{|U|^2} &  1+(U^0-1)\frac{(U^2)^2 }{|U|^2}&(U^0-1)\frac{U^2U^3}{|U|^2}  \cr
-U^3&  (U^0-1)\frac{U^1U^3}{|U|^2}& (U^0-1)\frac{U^2U^3}{|U|^2} &  1+(U^0-1)\frac{(U^3)^2}{|U|^2}
\end{bmatrix}
\end{align*}
which transforms $U^\mu=(\sqrt{c^2+|U|^2,U})$ into the local rest frame $(c,0,0,0).$ Applying the change of variables $P^\mu:=\Lambda p^\mu$, the left-hand side of \eqref{1} reduces to 
 \begin{align*}
\intr \inti f_E\Big(1+\frac{\mathcal{I}}{mc^2}\Big)^{-1} \phi(\I)\, d\I \frac{dp}{p^0}&=\frac{n}{M(\gamma)}\intr\inti  e^{-\frac{\gamma}{mc^2}\left(1+\frac{\mathcal{I}}{mc^2}\right)U^\mu p_\mu }\Big(1+\frac{\mathcal{I}}{mc^2}\Big)^{-1}\phi(\mathcal{I})  \,d\mathcal{I} \frac{dp}{p^0}\cr 
&=\frac{n}{M(\gamma)}\intr\inti  e^{-\frac{\gamma}{mc}\left(1+\frac{\mathcal{I}}{mc^2}\right)P^0 }\Big(1+\frac{\mathcal{I}}{mc^2}\Big)^{-1}\phi(\mathcal{I})  \,d\mathcal{I} \frac{dP}{P^0}
\end{align*}
where we used the fact that the volume element $\frac{dp}{p^0}$ and the Lorentz inner product $U^\mu p_\mu$ are invariant under the action of $\Lambda$. On the other hand, we observe that
\begin{equation}\label{P^0}
(mc)^2=p^\mu p_\mu=(\Lambda p^\mu)(\Lambda p_\mu)=(P^0)-|P|^2,\quad \mbox{and hence}\quad P^0=\sqrt{(mc)^2+|P|^2}.
\end{equation}
Thus one can rewrite $P$ to $ p$, and \eqref{1} becomes
$$
\frac{n}{M(\gamma)}\intr\inti  e^{-\frac{\gamma}{mc}\left(1+\frac{\mathcal{I}}{mc^2}\right)p^0 }\Big(1+\frac{\mathcal{I}}{mc^2}\Big)^{-1}\phi(\mathcal{I})  \,d\mathcal{I} \frac{dp}{p^0}=\intr\inti  f\Big(1+\frac{\I}{mc^2}\Big)^{-1}\phi(\mathcal{I}) \,d\mathcal{I} \frac{dp}{p^0}
$$
which gives the relation \eqref{gamma}. Next, to find $n$ and $U^\mu$, we apply the change of variables $P^\mu:=\Lambda p^\mu$ again to \eqref{2}. Then, in a similar way, the left-hand side reduces to 
\begin{align*}
\intr \inti p^\mu f_E \phi(\I)\, d\I \frac{dp}{p^0}&=\frac{n}{M(\gamma)}\intr\inti  p^\mu e^{-\frac{\gamma}{mc^2}\left(1+\frac{\mathcal{I}}{mc^2}\right)U^\mu p_\mu }\phi(\mathcal{I}) \,d\mathcal{I} \frac{dp}{p^0}\cr 
&= \frac{n}{M(\gamma)}\intr\inti  (\Lambda^{-1} P^\mu) e^{-\frac{\gamma}{mc}\left(1+\frac{\mathcal{I}}{mc^2}\right)P^0 }\phi(\mathcal{I}) \,d\mathcal{I} \frac{dP}{P^0}\cr 
&= \frac{n}{M(\gamma)}\Lambda^{-1}\intr\inti   p^\mu e^{-\frac{\gamma}{mc}\left(1+\frac{\mathcal{I}}{mc^2}\right)p^0 }\phi(\mathcal{I}) \,d\mathcal{I} \frac{dp}{p^0}.
\end{align*}
Since  $p^0=\sqrt{(mc)^2+|p|^2}$, it follows from the oddness  that
 
\begin{align*}
\intr\inti   p^\mu e^{-\frac{\gamma}{mc}\left(1+\frac{\mathcal{I}}{mc^2}\right)p^0 }\phi(\mathcal{I}) \,d\mathcal{I} \frac{dp}{p^0}&=\left(\intr\inti  e^{-\frac{\gamma}{mc}\left(1+\frac{\mathcal{I}}{mc^2}\right)p^0 }\phi(\mathcal{I}) \,d\mathcal{I} dp ,0,0,0\right)\cr 
&=\left(M(\gamma),0,0,0 \right),
\end{align*}
which gives
\begin{align} \label{(2)} \begin{split}
\frac{n}{M(\gamma)}\Lambda^{-1}\intr\inti   p^\mu e^{-\frac{\gamma}{mc}\left(1+\frac{\mathcal{I}}{mc^2}\right)p^0 }\phi(\mathcal{I}) \,d\mathcal{I} \frac{dp}{p^0}&= \frac{n}{M(\gamma)}\Lambda^{-1}\left(M(\gamma),0,0,0\right)\cr 
&= \frac nc\ \Lambda^{-1}(c,0,0,0)\cr 
&=\frac1c nU^\mu.
\end{split}\end{align}
In the last line, we used the fact that the Lorentz transformation $\Lambda: U^\mu\rightarrow (c,0,0,0)$ is invertible and thus $\Lambda^{-1}:(c,0,0,0)\rightarrow U^\mu$.  
On the other hand, the right-hand side of \eqref{2} can be expressed by \eqref{eckart} as
$$
\intr \inti p^\mu f\phi(\I)\, d\I \frac{dp}{p^0}=\frac 1c n_f U_f^\mu.
$$
Therefore, going back to \eqref{2} with \eqref{(2)}, we conclude that $n=n_f$ and $U^\mu=U_f^\mu$.
\noindent\newline
$\bullet$ Unique determination of $\gamma$: To prove this part, we will show the existence of a one-to-one correspondence between both sides of \eqref{gamma}. For this, it suffices to show that (1) $\widetilde{M}/M$ is strictly monotone in $\gamma \in(0,\infty)$, and (2) the ranges of both sides are the same. 
\noindent\newline
$(1)$ Strict monotonicity:  Using the change of variables $\frac{p}{mc}\rightarrow p$ and the spherical coordinates, one finds
\begin{align*}
\widetilde{M}(\gamma)&=4\pi(mc)^2\inti\inti  \frac{r^2}{\sqrt{1+r^2}} e^{-\left(1+\frac{\mathcal{I}}{mc^2}\right)\gamma\sqrt{1+r^2} }\Big(1+\frac{\I}{mc^2}\Big)^{-1} \phi(\mathcal{I})  \,d\mathcal{I} dr=: 4\pi (mc)^2 M_1(\gamma)
\end{align*}
and
\begin{align*}
M(\gamma)=4\pi (mc)^3\inti\inti  r^2 e^{-\left(1+\frac{\mathcal{I}}{mc^2}\right)\gamma\sqrt{1+r^2} } \phi(\mathcal{I})  \,d\mathcal{I} dr=:4\pi (mc)^3 M_2(\gamma).
\end{align*}
Also, it is straightforward that
\begin{align*}
\frac{d}{d\gamma} \{\widetilde{M}(\gamma) \}=-4\pi(mc)^2\inti\inti   r^2 e^{-\left(1+\frac{\mathcal{I}}{mc^2}\right)\gamma\sqrt{1+r^2} }\phi(\mathcal{I}) \,d\mathcal{I} dr =-4\pi (mc)^2M_2(\gamma)
\end{align*}
and
\begin{align*}
\frac{d}{d\gamma} \{M(\gamma) \}=-4\pi(mc)^3\inti\inti   r^2\sqrt{1+r^2} e^{-\left(1+\frac{\mathcal{I}}{mc^2}\right)\gamma\sqrt{1+r^2} }\Big(1+\frac{\mathcal{I}}{mc^2}\Big)\phi(\mathcal{I}) \,d\mathcal{I} dr=: -4\pi(mc)^3 M_3(\gamma).
\end{align*}
From these observations, we have  
\begin{align}\label{M/M} \begin{split}
\frac{d}{d\gamma}\biggl\{\frac{\widetilde{M}(\gamma)}{M(\gamma)}\biggl\}
&=\frac{\widetilde{M}^\prime(\gamma) M(\gamma)-\widetilde{M}(\gamma)M^\prime(\gamma)}{M^2(\gamma)}=\frac{1}{mc} \frac{M_1(\gamma)M_3(\gamma)-\{M_2(\gamma)\}^2}{\{M_2(\gamma)\}^2}
\end{split}\end{align}
where
\begin{align*}
&M_1(\gamma)M_3(\gamma)-\{M_2(\gamma)\}^2\cr 
&=\inti\inti  \frac{r^2}{\sqrt{1+r^2}} e^{-\left(1+\frac{\mathcal{I}}{mc^2}\right)\gamma\sqrt{1+r^2} }\frac{1}{1+\frac{\I}{mc^2}}\phi(\mathcal{I}) \,d\mathcal{I} dr\inti\inti   r^2\sqrt{1+r^2} e^{-\left(1+\frac{\mathcal{I}}{mc^2}\right)\gamma\sqrt{1+r^2} }\phi(\mathcal{I})\left(1+\frac{\mathcal{I}}{mc^2}\right) \,d\mathcal{I} dr\cr 
&-\left(\inti\inti   r^2 e^{-\left(1+\frac{\mathcal{I}}{mc^2}\right)\gamma\sqrt{1+r^2} }\phi(\mathcal{I}) \,d\mathcal{I} dr\right)^2.
\end{align*}
Therefore, since the integrands above are linearly independent, by H\"{o}lder's inequality we conclude that $\widetilde{M}/M$  is strictly increasing in $\gamma\in (0,\infty)$.
\noindent\newline
$(2)$ Ranges: It follows from \eqref{n} and the change of variables $P^\beta:=\Lambda p^\beta$ that
\begin{align*}
\frac {1}{n_f}\intr\inti  f\Big(1+\frac{\I}{mc^2}\Big)^{-1}\phi(\mathcal{I}) \,d\mathcal{I} \frac{dp}{p^0}&=c\frac {\intr\inti  f\Big(1+\frac{\I}{mc^2}\Big)^{-1}\phi(\mathcal{I}) \,d\mathcal{I} \frac{dp}{p^0}}{ \int_{\mathbb{R}^3}\int_0^\infty U_{f\mu} p^\mu f \phi(\mathcal{I}) \,d\mathcal{I}\,\frac{dp}{p^0}}\cr  
&=\frac {\intr\inti  f_\Lambda\Big(1+\frac{\I}{mc^2}\Big)^{-1}\phi(\mathcal{I}) \,d\mathcal{I} \frac{dP}{P^0}}{ \int_{\mathbb{R}^3}\int_0^\infty  f_\Lambda \phi(\mathcal{I}) \,d\mathcal{I}\,dP}\cr 
&=\frac{1}{mc}\frac {\intr\inti  \frac{1}{\sqrt{1+|\frac{P}{mc}|^2}}f_\Lambda\Big(1+\frac{\I}{mc^2}\Big)^{-1}\phi(\mathcal{I}) \,d\mathcal{I} \,dP}{ \int_{\mathbb{R}^3}\int_0^\infty  f_\Lambda \phi(\mathcal{I}) \,d\mathcal{I}\,dP}
\end{align*}
where $f_\Lambda:=f(x^\alpha,\Lambda^{-1}P^\beta)$ and  we used \eqref{P^0} in the last line. Since the above is strictly less than $\frac{1}{mc}$, it only remains to show that $\mbox{Range}(\widetilde{M}/M)=(0,\frac{1}{mc})$. For this, we employ the integration by parts with respect to $r$ to see
\begin{align*}
\inti\inti   r^2 e^{-\left(1+\frac{\mathcal{I}}{mc^2}\right)\gamma\sqrt{1+r^2} }\phi(\mathcal{I}) \,d\mathcal{I} dr&=\frac{1}{\gamma}\inti\inti \left\{ \sqrt{1+r^2}+\frac{r^2}{\sqrt{1+r^2}} \right\}   e^{-\left(1+\frac{\mathcal{I}}{mc^2}\right)\gamma\sqrt{1+r^2} }\frac{1}{1+\frac{\I}{mc^2}}\phi(\mathcal{I}) \,d\mathcal{I} dr\cr 
&\ge \frac{1}{\gamma}\inti\inti  \frac{r^2}{\sqrt{1+r^2}} e^{-\left(1+\frac{\mathcal{I}}{mc^2}\right)\gamma\sqrt{1+r^2} }\frac{1}{1+\frac{\I}{mc^2}}\phi(\mathcal{I}) \,d\mathcal{I} dr, 
\end{align*}
which combined with \eqref{M/M} gives
$$
\frac{\widetilde{M}(\gamma)}{M(\gamma)} \le \frac{\gamma}{mc},\quad\mbox{and hence}\quad \frac{\widetilde{M}(\gamma)}{M(\gamma)} \rightarrow 0 \ \ \mbox{as}\ \ \gamma\rightarrow 0.
$$
On the other hand, it follows from the H\"{o}lder inequality that
\begin{align*}
\frac{\widetilde{M}(\gamma)}{M(\gamma)}&=\frac{1}{mc}\frac{\inti\inti  \frac{r^2}{\sqrt{1+r^2}} e^{-\left(1+\frac{\mathcal{I}}{mc^2}\right)\gamma\sqrt{1+r^2} }\frac{1}{1+\frac{\I}{mc^2}}\phi(\mathcal{I}) \,d\mathcal{I} dr}{\inti\inti   r^2 e^{-\left(1+\frac{\mathcal{I}}{mc^2}\right)\gamma\sqrt{1+r^2} }\phi(\mathcal{I}) \,d\mathcal{I} dr}\cr 
&\ge \frac{1}{mc}\left(\frac{\inti\inti  \frac{r^2}{\sqrt{1+r^2}} e^{-\left(1+\frac{\mathcal{I}}{mc^2}\right)\gamma\sqrt{1+r^2} }\frac{1}{1+\frac{\I}{mc^2}}\phi(\mathcal{I}) \,d\mathcal{I} dr}{\inti\inti   r^2\sqrt{1+r^2} e^{-\left(1+\frac{\mathcal{I}}{mc^2}\right)\gamma\sqrt{1+r^2} }(1+\frac{\I}{mc^2}) \phi(\mathcal{I}) \,d\mathcal{I} dr}\right)^{1/2}
\end{align*}
Since $\frac{1}{1+x}\ge  (1+x)-2x$ for all $x\ge 0,$ we obtain
\begin{align}\label{I+II}\begin{split}
\frac{\widetilde{M}(\gamma)}{M(\gamma)}&\ge \frac{1}{mc} \Biggl(\frac{\inti\inti  \frac{r^2}{\sqrt{1+r^2}} e^{-\left(1+\frac{\mathcal{I}}{mc^2}\right)\gamma\sqrt{1+r^2} }(1+\frac{\I}{mc^2}) \phi(\mathcal{I}) \,d\mathcal{I} dr}{\inti\inti   r^2\sqrt{1+r^2} e^{-\left(1+\frac{\mathcal{I}}{mc^2}\right)\gamma\sqrt{1+r^2} }(1+\frac{\I}{mc^2}) \phi(\mathcal{I}) \,d\mathcal{I} dr}\cr 
&\qquad -\frac{2}{mc^2}\frac{\inti\inti  \frac{r^2}{\sqrt{1+r^2}} e^{-\left(1+\frac{\mathcal{I}}{mc^2}\right)\gamma\sqrt{1+r^2} } \I\phi(\mathcal{I})   \,d\mathcal{I} dr}{\inti\inti   r^2\sqrt{1+r^2} e^{-\left(1+\frac{\mathcal{I}}{mc^2}\right)\gamma\sqrt{1+r^2} }(1+\frac{\I}{mc^2}) \phi(\mathcal{I}) \,d\mathcal{I} dr}\Biggl)^{1/2}\cr 
&=:\frac {1}{mc}(I-II)^{\frac 12}.
\end{split}\end{align}
For $I$, we use the identity
$$
\frac{1}{\sqrt{1+r^2}}=\sqrt{1+r^2}-\frac{r^2}{\sqrt{1+r^2}}
$$ 
to obtain
\begin{align*}
I&=1-\frac{\inti\inti  \frac{r^4}{\sqrt{1+r^2}} e^{-\left(1+\frac{\mathcal{I}}{mc^2}\right)\gamma\sqrt{1+r^2} }(1+\frac{\I}{mc^2}) \phi(\mathcal{I}) \,d\mathcal{I} dr}{\inti\inti   r^2\sqrt{1+r^2} e^{-\left(1+\frac{\mathcal{I}}{mc^2}\right)\gamma\sqrt{1+r^2} }(1+\frac{\I}{mc^2}) \phi(\mathcal{I}) \,d\mathcal{I} dr}.
\end{align*}
Applying the integration by parts with respect to $r$,  we get
\begin{align*}
I&=1-\frac{3}{\gamma}\frac{\inti\inti  r^2 e^{-\left(1+\frac{\mathcal{I}}{mc^2}\right)\gamma\sqrt{1+r^2} }\phi(\mathcal{I}) \,d\mathcal{I} dr}{\inti\inti   r^2\sqrt{1+r^2} e^{-\left(1+\frac{\mathcal{I}}{mc^2}\right)\gamma\sqrt{1+r^2} }(1+\frac{\I}{mc^2}) \phi(\mathcal{I}) \,d\mathcal{I} dr}\cr 
&\ge 1-\frac 3\gamma.
\end{align*}
For $II$,  since we have set $\phi(\I)=\I^\sigma$ $(\sigma> -1)$, applying the integration by parts with respect to $\I$ gives
\begin{align}\label{phi lower} \begin{split}
II&= \frac{2}{mc^2}\frac{\inti\inti  \frac{r^2}{\sqrt{1+r^2}} e^{-\left(1+\frac{\mathcal{I}}{mc^2}\right)\gamma\sqrt{1+r^2} }\phi(\mathcal{I})\I  \,d\mathcal{I} dr}{\inti\inti   r^2\sqrt{1+r^2} e^{-\left(1+\frac{\mathcal{I}}{mc^2}\right)\gamma\sqrt{1+r^2} }(1+\frac{\I}{mc^2}) \phi(\mathcal{I}) \,d\mathcal{I} dr} \cr
&= \frac{2(\sigma+1)}{\gamma}\frac{\inti\inti  \frac{r^2}{ 1+r^2} e^{-\left(1+\frac{\mathcal{I}}{mc^2}\right)\gamma\sqrt{1+r^2} }\phi(\mathcal{I})  \,d\mathcal{I} dr}{\inti\inti   r^2\sqrt{1+r^2} e^{-\left(1+\frac{\mathcal{I}}{mc^2}\right)\gamma\sqrt{1+r^2} }(1+\frac{\I}{mc^2}) \phi(\mathcal{I}) \,d\mathcal{I} dr}\cr 
&\le  \frac{2(\sigma+1)}{\gamma}.
\end{split}\end{align}
Going back to \eqref{I+II} with estimates of $I$ and $II$, we conclude that
$$
\frac{\widetilde{M}(\gamma)}{M(\gamma)}\ge \frac{1}{mc} \left( 1-\frac3\gamma - \frac{2(\sigma+1)}{\gamma}\right)^{\frac 12},\qquad\mbox{and hence}\qquad \frac{\widetilde{M}(\gamma)}{M(\gamma)}\rightarrow \frac{1}{mc} \ \ \mbox{as}\ \  \gamma\rightarrow \infty.
$$
This completes the proof.

%
%
%
%
%
%



	%
	%
	%
	%

\end{document}